\documentclass[12pt,english]{amsart}
\usepackage{mathptmx}
\usepackage[T1]{fontenc}
\usepackage[latin9]{inputenc}
\usepackage[a4paper]{geometry}
\usepackage{prettyref}
\usepackage{amsmath}
\usepackage{amssymb}

\makeatletter

\newcommand{\noun}[1]{\textsc{#1}}

\usepackage{latexsym, pifont, color, graphicx, mflogo}
\theoremstyle{plain}
\newtheorem{stthm}{Theorem}[section]

\numberwithin{equation}{section} 
\numberwithin{figure}{section} 
\numberwithin{table}{section} 

\def\newrefformat#1#2{%
  \@namedef{pr@#1}##1{#2}}
\newrefformat{eq}{\textup{(\ref{#1})}}
\newrefformat{ssec}{sub-section \ref{#1}}
\newrefformat{sssec}{sub-sub-section \ref{#1}}
\newrefformat{cha}{chapter \ref{#1}}
\newrefformat{sec}{section \ref{#1}}
\newrefformat{tab}{table \ref{#1} on page \pageref{#1}}
\newrefformat{fig}{figure \ref{#1} on page \pageref{#1}}
\def\prettyref#1{\@prettyref#1|}
\def\@prettyref#1|#2|{%
  \expandafter\ifx\csname pr@#1\endcsname\relax%
    \PackageWarning{prettyref}{Reference format #1\space undefined}%
    \ref{#1|#2}%
  \else%
    \csname pr@#1\endcsname{#1|#2}%
  \fi%
}

\DeclareRobustCommand{\QED}{%
  \ifmmode 
  \else \leavevmode\unskip\penalty9999 \hbox{}\nobreak\hfill
  \fi
  \quad\hbox{\QEDsymbol}}
\newcommand{\QEDsymbol}{Q.E.D.}

 \theoremstyle{remark}
 \newtheorem*{notation}{Notation}
 \theoremstyle{plain}
 \newtheorem{stprop}[stthm]{Proposition}
 \theoremstyle{plain}
 \newtheorem{stlem}[stthm]{Lemma}
 \theoremstyle{remark}
 \newtheorem{strem}[stthm]{Remark}
 \theoremstyle{plain}
 \newtheorem{stcor}[stthm]{Corollary}

\usepackage[bookmarksopen, naturalnames]{hyperref}
\newrefformat{thm}{Theorem \ref{#1}}
\newrefformat{prop}{Proposition \ref{#1}}
\newrefformat{cor}{Corollary \ref{#1}}
\newrefformat{lem}{Lemma \ref{#1}}
\newrefformat{rem}{Remark \ref{#1}}
\newrefformat{sec}{Section \ref{#1}}
\newrefformat{sub}{Subsection \ref{#1}}
\newrefformat{subsub}{Paragraph \ref{#1}}
\newrefformat{def}{Definition \ref{#1}}

{\end{list}}

\makeatother

\usepackage{babel}

\begin{document}

\title{Compact composition operators on Hardy-Orlicz and weighted Bergman-Orlicz
spaces on the ball}

\author{Stéphane Charpentier}
\begin{abstract}
Using recent characterizations of the compactness of composition operators
on Hardy-Orlicz and Bergman-Orlicz spaces on the ball (\cite{Charp1,Charp2}),
we first show that a composition operator which is compact on every
Hardy-Orlicz (or Bergman-Orlicz) space has to be compact on $H^{\infty}$.
Then, although it is well-known that a map whose range is contained
in some nice Kor\'anyi approach region induces a compact composition
operator on $H^{p}\left(\mathbb{B}_{N}\right)$ or on $A_{\alpha}^{p}\left(\mathbb{B}_{N}\right)$,
we prove that, for each Kor\'anyi region $\Gamma$, there exists
a map $\phi:\mathbb{B}_{N}\rightarrow\Gamma$ such that, $C_{\phi}$
is not compact on $H^{\psi}\left(\mathbb{B}_{N}\right)$, when $\psi$
grows fast. Finally, we extend (and simplify the proof of) a result
by K. Zhu for classical weighted Bergman spaces, by showing that,
under reasonable conditions, a composition operator $C_{\phi}$ is
compact on the weighted Bergman-Orlicz space $A_{\alpha}^{\psi}\left(\mathbb{B}_{N}\right)$,
if and only if\[
\lim_{\left|z\right|\rightarrow1}\frac{\psi^{-1}\left(1/\left(1-\left|\phi\left(z\right)\right|\right)^{N\left(\alpha\right)}\right)}{\psi^{-1}\left(1/\left(1-\left|z\right|\right)^{N\left(\alpha\right)}\right)}=0.\]
In particular, we deduce that the compactness of composition operators
on $A_{\alpha}^{\psi}\left(\mathbb{B}_{N}\right)$ does not depend
on $\alpha$ anymore when the Orlicz function $\psi$ grows fast.
\end{abstract}



\address{Charpentier Stéphane, Département de Mathématiques, Bâtiment 425,
Université Paris-Sud, F-91405, Orsay, France}

\email{stephane.charpentier@math.u-psud.fr}

\maketitle
\def\leftmark{\MakeUppercase{Compact composition operators on Hardy-Orlicz
and Bergman-Orlicz spaces on $\mathbb{B}_{N}$}}

\def\rightmark{\MakeUppercase{Stéphane Charpentier}}

\section{Introduction}

Let $\mathbb{B}_{N}=\left\{ z=\left(z_{1},\ldots z_{N}\right)\in\mathbb{C}^{N},\,\sum_{i=1}^{N}\left|z_{i}\right|^{2}<1\right\} $
denote the open unit ball of $\mathbb{C}^{N}$. Given a holomorphic
map $\phi:\mathbb{B}_{N}\rightarrow\mathbb{B}_{N}$, the composition
operator $C_{\phi}$ of symbol $\phi$ is defined by $C_{\phi}\left(f\right)=f\circ\phi$,
for $f$ holomorphic on $\mathbb{B}_{N}$. Composition operators have
been extensively studied on common Banach spaces of analytic functions,
in particular on the Hardy spaces $H^{p}\left(\mathbb{B}_{N}\right)$
and on the Bergman spaces $A^{p}\left(\mathbb{B}_{N}\right)$, $1\leq p<\infty$.
The continuity and compactness of these operators have been characterized
in terms of Carleson measures (\cite{COWEN-MACCLUER}). In dimension
one, the boundedness of $C_{\phi}$ for any $\phi:\mathbb{D}\rightarrow\mathbb{D}$
is a consequence of the Littlewood subordination principle (\cite{SHAPIRO}).
In $\mathbb{C}^{N}$, $N>1$, it is well-known that there exists some
map $\phi:\mathbb{B}_{N}\rightarrow\mathbb{B}_{N}$ such that the
associated composition operator is not bounded on $H^{p}\left(\mathbb{B}_{N}\right)$.
Whatever the dimension, it appears that both boundedness and compactness
of $C_{\phi}$ on $H^{p}\left(\mathbb{B}_{N}\right)$ (resp. $A^{p}\left(\mathbb{B}_{N}\right)$)
are independent of $p$. On the other hand, every composition operator
is obviously bounded on $H^{\infty}$ and it is not difficult to check
that $C_{\phi}$ is compact on $H^{\infty}$ if and only if $\left\Vert \phi\right\Vert _{\infty}<1$.
Thus there is a {}``break'' between $H^{\infty}$ and $H^{p}\left(\mathbb{B}_{N}\right)$
(resp. $A^{p}\left(\mathbb{B}_{N}\right)$), for the compactness in
dimension one, and even for the boundedness, when $N>1$.

These observations first motivated P. Lefèvre, D. Li, H. Queffélec
and L. Rodr\'iguez-Piazza to study composition operators on Hardy-Orlicz
spaces $H^{\psi}\left(\mathbb{D}\right)$ (resp. Bergman-Orlicz spaces
$A^{\psi}\left(\mathbb{D}\right)$) of the disc (\cite{QUEF-LI-LE-RO-PI,QUEF-LI-LEF-RO-B-O-H-O,QUEF-LI-LEF-ROD-H-2-H-O-COMP-OP,QUEF-LI-LEFEVRE-BERGMAN-ORLICZ}),
and then the author of \cite{Charp1,Charp2} to look at these questions
in $\mathbb{C}^{N}$. These spaces both provide an intermediate scale
of spaces between $H^{\infty}$ and $H^{p}\left(\mathbb{B}_{N}\right)$
(resp. $A^{p}\left(\mathbb{B}_{N}\right)$) and generalize the latter.
In particular, in \cite{QUEF-LI-LEF-RO-B-O-H-O}, the authors were
interested in the question of whether there are some Hardy-Orlicz
spaces on which the compactness of $C_{\phi}$ is equivalent to that
on $H^{\infty}$. In fact, they answer this question in the negative,
by proving (\cite[Theorem 4.1]{QUEF-LI-LEF-RO-B-O-H-O}) that, for
every Hardy-Orlicz space $H^{\psi}\left(\mathbb{D}\right)$, one can
construct a surjective map $\phi:\mathbb{D}\rightarrow\mathbb{D}$
which induces a compact composition operator $C_{\phi}$ on $H^{\psi}\left(\mathbb{D}\right)$.
This result extends that obtained by B. MacCluer and J. Shapiro for
$H^{p}\left(\mathbb{D}\right)$ (\cite[Example 3.12]{MACCLUER-SHAPIRO}).
The same problem in the Bergman-Orlicz case has not yet been completely
solved. In several variables, the situation is much more surprizing,
as we show in \cite{Charp1,Charp2} that there exist some Hardy-Orlicz
and Bergman-Orlicz spaces, {}``close'' enough to $H^{\infty}$,
on which every composition operator is bounded.

In this paper, we are mainly interested in the possibility to extend
some known results about compactness of composition operators on classical
Hardy or Bergman spaces, to the corresponding Orlicz spaces. We think
that this study may outline some interesting phenomena and precise
the link between the behavior of $C_{\phi}$ and that of $\phi$.

First of all, we come back to the {}``break'' between $H^{\infty}$
and $H^{p}$, $1\leq p<\infty$, for the compactness of $C_{\phi}$.
There is no difference between being compact for $C_{\phi}$ on one
$H^{p}\left(\mathbb{B}_{N}\right)$ and on every $H^{p}\left(\mathbb{B}_{N}\right)$,
while this property clearly depends on the Orlicz function $\psi$
in $H^{\psi}\left(\mathbb{B}_{N}\right)$. Therefore, we can wonder
if the above question answered by \cite{QUEF-LI-LEF-RO-B-O-H-O} was
the good one; indeed, the study of $C_{\phi}$ on Hardy-Orlicz spaces
arises the following question: what can we say about a composition
operator which is compact on every Hardy-Orlicz space? It turns out
that such an operator has to be compact on $H^{\infty}$, which seems
to us to be a positive result, because it confirms that Hardy-Orlicz
spaces covers well the {}``gap'' between every $H^{p}$ and $H^{\infty}$.
This result also stands when we replace Hardy-Orlicz spaces by Bergman-Orlicz
spaces.

Moreover, in Hardy or Bergman spaces, compactness (and boundedness)
of composition operators is handled in terms of geometric conditions,
emphasizing the importance of the manner in which the symbol $\phi$
approaches the boundary of $\mathbb{B}_{N}$. To be precise, if we
denote by $\Gamma\left(\zeta,a\right)\subset\mathbb{B}_{N}$, for
$\zeta\in\mathbb{S}_{N}$ and $a>1$, the Kor\'anyi approach region\[
\Gamma\left(\zeta,a\right)=\left\{ z\in\mathbb{B}_{N},\,\left|1-\left\langle z,\zeta\right\rangle \right|<\frac{a}{2}\left(1-\left|z\right|^{2}\right)\right\} ,\]
it is known (\cite{MacCluer_95}) that if $\phi$ takes the unit ball
into a Kor\'anyi region $\Gamma\left(\zeta,a\right)$ with a small
enough angular opening $a$, then $C_{\phi}$ is compact on $H^{p}\left(\mathbb{B}_{N}\right)$
and on $A_{\alpha}^{p}\left(\mathbb{B}_{N}\right)$. When $N=1$,
the Kor\'anyi regions are just non-tangential approach regions. In
this paper, we show that this result does not hold anymore for Hardy-Orlicz
spaces on $\mathbb{B}_{N}$; for Bergman-Orlicz spaces, we obtain
such a result in dimension one only.

In \cite{MACCLUER-SHAPIRO}, the authors related the compactness of
the composition operator $C_{\phi}$ on $H^{p}\left(\mathbb{D}\right)$
or $A_{\alpha}^{p}\left(\mathbb{D}\right)$ to the existence of angular
derivative for $\phi$ at the boundary. We say that the \emph{angular
derivative} of $\phi$ exists at a point $\zeta\in\mathbb{T}$ if
there exists $\omega\in\mathbb{T}$ such that\[
\frac{\phi\left(z\right)-\omega}{z-\zeta}\]
has a finite limit as $\zeta$ tends non-tangentially to $\zeta$
through $\mathbb{D}$. The Julia-Caratheodory Theorem then asserts
that the non-existence of an angular derivative for $\phi$ at some
$\zeta\in\mathbb{T}$ is equivalent to\begin{equation}
\lim_{z\rightarrow\zeta}\frac{1-\left|z\right|}{1-\left|\phi\left(z\right)\right|}=0.\label{eq|ang_der_J_C_Thm}\end{equation}
Shapiro and Taylor \cite{SHAPIRO-TAYLOR} pointed out that if $C_{\phi}$
is to be compact on $H^{p}\left(\mathbb{D}\right)$, then $\phi$
cannot have an angular derivative at even a single point in $\mathbb{T}$,
which may be written:\begin{equation}
\lim_{\left|z\right|\rightarrow1}\frac{1-\left|z\right|}{1-\left|\phi\left(z\right)\right|}=0.\label{eq|eq_carac_Zhu_intro}\end{equation}
In \cite{MACCLUER-SHAPIRO}, it is proven that \prettyref{eq|ang_der_J_C_Thm}
is not sufficient to the compactness of $C_{\phi}$ on Hardy spaces
of the unit disc in general, yet it is when $\phi$ is univalent.
However, this condition is necessary and sufficient for $C_{\phi}$
to be compact on every weighted Bergman spaces of the disc. The last
main goal of this paper is to extend some of these results to Hardy-Orlicz
and Bergman-Orlicz spaces of the unit ball.

In several variables, we can also define the angular derivative of
$\phi:\mathbb{B}_{N}\rightarrow\mathbb{B}_{N}$ at a point in the
unit sphere $\mathbb{S}_{N}$ and the Julia-Caratheodory Theorem also
holds in $\mathbb{B}_{N}$ (\cite[Theorem 8.5.6]{RUDIN_80}). Here,
as we already said, the situation is complicated by the fact that
some composition operators are not bounded on Hardy or Bergman spaces,
and the fact that even the boundedness of $C_{\phi}$ on $A_{\alpha}^{p}\left(\mathbb{B}_{N}\right)$
depends on $\alpha$. In \cite{ZHU_COMP_OP_BERG_BALL}, K. Zhu proves
that $C_{\phi}$ is compact on $A_{\alpha}^{p}\left(\mathbb{B}_{N}\right)$
if and only if Condition \prettyref{eq|eq_carac_Zhu_intro} is satisfied,
whenever $C_{\phi}$ is bounded on some $A_{\beta}^{p}\left(\mathbb{B}_{N}\right)$,
for some $-1<\beta<\alpha$. This assumption is somehow justified
by the above observation and by \cite[Section 6]{MACCLUER-SHAPIRO},
in which the authors show that, for any $\alpha>-1$ and any $0<p<\infty$,
there exists $\phi:\mathbb{B}_{N}\rightarrow\mathbb{B}_{N}$ with
no angular derivative at any point of $\mathbb{S}_{N}$, such that
$C_{\phi}$ is bounded on $A_{\alpha}^{p}\left(\mathbb{B}_{N}\right)$
but not compact. There even exists such a map $\phi$ such that $C_{\phi}$
is not bounded on $A_{\alpha}^{p}\left(\mathbb{B}_{N}\right)$. In
the present paper, we generalize Zhu's result to weighted Bergman-Orlicz
spaces on the ball, by using recent characterizations of boundedness
and compactness of composition operators on these spaces (\cite{Charp1}).
We show that, if $C_{\phi}$ is bounded on some $A_{\beta}^{\psi}\left(\mathbb{B}_{N}\right)$,
$-1<\beta<\alpha$, then it is compact on $A_{\alpha}^{\psi}\left(\mathbb{B}_{N}\right)$
if and only if\begin{equation}
\lim_{\left|z\right|\rightarrow1}\frac{\psi^{-1}\left(1/\left(1-\left|\phi\left(z\right)\right|\right)^{N\left(\alpha\right)}\right)}{\psi^{-1}\left(1/\left(1-\left|z\right|\right)^{N\left(\alpha\right)}\right)}=0,\label{eq|eq_MAIN_3_INTRO}\end{equation}
where $N\left(\alpha\right)=N+\alpha+1$, under a mild and usual regularity
condition on the Orlicz function $\psi$. Our proof is quite simple,
while that of K. Zhu uses a Schur test in $H^{2}\left(\mathbb{B}_{N}\right)$
and the fact that the compactness of composition operators on $H^{p}\left(\mathbb{B}_{N}\right)$
does not depend on $p$. Combining this result with the automatic
boundedness of every composition operator on $A_{\alpha}^{\psi}\left(\mathbb{B}_{N}\right)$
when $\psi$ satisfies the $\Delta^{2}$-Condition, we get that the
compactness on such $A_{\alpha}^{\psi}\left(\mathbb{B}_{N}\right)$
does not depend on $\alpha$ anymore. To be precise, $C_{\phi}$ is
compact on $A_{\alpha}^{\psi}\left(\mathbb{B}_{N}\right)$ if and
only if\[
\lim_{\left|z\right|\rightarrow1}\frac{\psi^{-1}\left(1/\left(1-\left|\phi\left(z\right)\right|\right)\right)}{\psi^{-1}\left(1/\left(1-\left|z\right|\right)\right)}=0,\]
whenever $\psi$ satisfies the $\Delta^{2}$-Condition.

We have to mention that Condition \prettyref{eq|eq_MAIN_3_INTRO}
is, in any case, necessary. Moreover, the authors of \cite{QUEF-LI-LE-RO-PI}
obtained such a result in dimension one, as announced in \cite{LI}.
However, their proof uses the characterization of the compactness
of composition operators in terms of the Nevanlinna counting function
and is more complicated.

We organize our paper as follows: a first preliminary part is devoted
to the definitions and the statements of the already known results
we need. The main part contains the three most important results mentionned
above.

\smallskip{}

\begin{notation}
Throughout this paper, we will denote by $d\sigma_{N}$ the normalized
invariant measure on the unit sphere $\mathbb{S}_{N}=\partial\mathbb{B}_{N}$,
and by $dv_{\alpha}=c_{\alpha}\left(1-\left|z\right|^{2}\right)^{\alpha}dv$,
$\alpha>-1$, the normalized weighted Lebesgue measure on the ball.

Given two points $z,w\in\mathbb{C}^{N}$, the euclidean inner product
of $z$ and $w$ will be denoted by $\left\langle z,w\right\rangle $,
that is $\left\langle z,w\right\rangle =\sum_{i=1}^{N}z_{i}\overline{w_{i}}$;
the notation $\left|\cdot\right|$ will stand for the associated norm,
as well as for the modulus of a complex number.

If $\alpha>-1$ is a real number, we will denote by $N\left(\alpha\right)$
the quantity $N+\alpha+1$.
\end{notation}

\section{\label{sub|preliminaries}Preliminaries}

\subsection{Hardy-Orlicz and Bergman-Orlicz spaces - Definitions}

A strictly convex function $\psi:\mathbb{R}_{+}\rightarrow\mathbb{R}_{+}$
is called an Orlicz function if $\psi\left(0\right)=0$, $\psi$ is
continuous at $0$ and ${\displaystyle \frac{\psi\left(x\right)}{x}\xrightarrow[x\rightarrow+\infty]{}+\infty}$.
If $\left(\Omega,\mathbb{P}\right)$ is a probability space, the Orlicz
space $L^{\psi}\left(\Omega\right)$ associated to the Orlicz function
$\psi$ on $\left(\Omega,\mathbb{P}\right)$ is the set of all (equivalence
classes of) measurable functions $f$ on $\Omega$ such that there
exists some $C>0$, such that $\int_{\Omega}\psi\left(\frac{\left|f\right|}{C}\right)d\mathbb{P}$
is finite. $L^{\psi}\left(\Omega\right)$ is a vector space, which
can be normed with the so-called Luxemburg norm defined by\[
\left\Vert f\right\Vert _{\psi}=\inf\left\{ C>0,\,\int_{\Omega}\psi\left(\frac{\left|f\right|}{C}\right)d\mathbb{P}\leq1\right\} .\]
It is well-known that $\left(L^{\psi}\left(\Omega\right),\left\Vert \cdot\right\Vert _{\psi}\right)$
is a Banach space (see \cite{RAO-REN}).

\smallskip{}
Taking $\Omega=\mathbb{S}_{N}$ and $d\mathbb{P}=d\sigma_{N}$, the
Hardy-Orlicz space $H^{\psi}\left(\mathbb{B}_{N}\right)$ on $\mathbb{B}_{N}$
is the Banach space of analytic functions $f:\mathbb{B}_{N}\rightarrow\mathbb{C}$
such that $\left\Vert f\right\Vert _{H^{\psi}}:=\sup_{0<r<1}\left\Vert f_{r}\right\Vert _{\psi}<\infty$,
where $f_{r}\in L^{\psi}\left(\mathbb{S}_{N}\right)$ is defined by
$f_{r}\left(z\right)=f\left(rz\right)$. Every function $f\in H^{\psi}\left(\mathbb{B}_{N}\right)$
admits a radial boundary limit $f^{*}$ such that $\left\Vert f^{*}\right\Vert _{\psi}=\sup_{0<r<1}\left\Vert f_{r}\right\Vert _{\psi}<\infty$
(\cite[Section 1.3]{Charp2}). For simplicity, we will sometimes denote
by $\|\cdot\|_{\psi}$ the norm on $H^{\psi}\left(\mathbb{B}_{N}\right)$,
emphasizing that $H^{\psi}\left(\mathbb{B}_{N}\right)$ can be seen
as a subspace of $L^{\psi}\left(\mathbb{S}_{N}\right)$.

\smallskip{}
With $\Omega=\mathbb{B}_{N}$ and $d\mathbb{P}=dv_{\alpha}$, $\alpha>-1$,
the weighted Bergman-Orlicz space $A_{\alpha}^{\psi}\left(\mathbb{B}_{N}\right)$
is $L^{\psi}\left(\mathbb{B}_{N}\right)\cap H\left(\mathbb{B}_{N}\right)$,
where $H\left(\mathbb{B}_{N}\right)$ stands for the vector space
of analytic functions on the unit ball. $A_{\alpha}^{\psi}\left(\mathbb{B}_{N}\right)$
is a Banach space.

\smallskip{}
From the definitions, it is easy to verify that the following inclusions
hold:\[
H^{\infty}\subset H^{\psi}\left(\mathbb{B}_{N}\right)\subset H^{p}\left(\mathbb{B}_{N}\right)\mbox{ and }H^{\infty}\subset A_{\alpha}^{\psi}\left(\mathbb{B}_{N}\right)\subset A_{\alpha}^{p}\left(\mathbb{B}_{N}\right)\]
for every Orlicz function $\psi$ and any $1\leq p<\infty$. Moreover,
if $\psi\left(x\right)=x^{p}$, for some $1\leq p<\infty$ and for
every $x\geq0$, then $H^{\psi}\left(\mathbb{B}_{N}\right)=H^{p}\left(\mathbb{B}_{N}\right)$
and $A_{\alpha}^{\psi}\left(\mathbb{B}_{N}\right)=A_{\alpha}^{p}\left(\mathbb{B}_{N}\right)$.

\smallskip{}

\subsection{Four classes of Orlicz functions}

Let $\psi$ be an Orlicz function. In order to distinguish the Orlicz
spaces and to get a significant scale of intermediate spaces between
$L^{\infty}$ and $L^{p}\left(\Omega\right)$, we define four classes
of Orlicz functions.

\smallskip{}

\---\enskip The two first conditions are regularity conditions: we
say that $\psi$ satisfies the $\nabla_{0}$-condition if it satisfies
one of the following two equivalent conditions:\begin{enumerate}\renewcommand{\theenumi}{\roman{enumi}}\item For
any $B>1$, there exists some constant $C_{B}\geq1$, such that ${\displaystyle \frac{\psi\left(Bx\right)}{\psi\left(x\right)}\leq\frac{\psi\left(C_{B}By\right)}{\psi\left(y\right)}}$
for any $x\leq y$ large enough;\item For any $n>0$, there exists
$C_{n}>0$ such that ${\displaystyle \frac{\psi\left(Bx\right)^{n}}{\psi\left(x\right)^{n}}\leq\frac{\psi\left(C_{n}y\right)}{\psi\left(By\right)}}$
for any $x\leq y$ large enough.\end{enumerate}

Let us notice that (2)$\,\Rightarrow\,$(1) is obvious, while an easy
induction allows to prove (1)$\,\Rightarrow\,$(2); the details are
left to the reader.

If the constant $C_{B}$ can be chosen independently of $B$, then
$\psi$ satisfies the \emph{uniform} $\nabla_{0}$-Condition.\medskip{}

\---\enskip The $\nabla_{2}$-class consists of those Orlicz functions
$\psi$ such that there exist some $\beta>1$ and some $x_{0}>0$,
such that $\psi\left(\beta x\right)\geq2\beta\psi\left(x\right)$,
for $x\geq x_{0}$.\medskip{}

\---\enskip The third one is a condition of moderate growth: $\psi$
satisfies the $\Delta_{2}$-Condition if there exist $x_{0}>0$ and
a constant $K>1$, such that $\psi\left(2x\right)\leq K\psi\left(x\right)$
for any $x\geq x_{0}$.

\medskip{}

\---\enskip The fourth condition is a condition of fast growth: $\psi$
satisfies the $\Delta^{2}$-Condition if it satisfies one of the following
equivalent conditions:\begin{enumerate}\renewcommand{\theenumi}{\roman{enumi}}\item There
exist $C>0$ and $x_{0}>0$, such that $\psi\left(x\right)^{2}\leq\psi\left(Cx\right)$
for every $x\geq x_{0}$;\item There exist $b>1$, $C>0$ and $x_{0}>0$
such that $\psi\left(x\right)^{b}\leq\psi\left(Cx\right)$, for every
$x\geq x_{0}$;\item For every $b>1$, there exist $C_{b}>0$ and
$x_{0,b}>0$ such that $\psi\left(x\right)^{b}\leq\psi\left(C_{b}x\right)$,
for every $x\geq x_{0,b}$.\end{enumerate}

\medskip{}

Finally, we mention that these conditions are not independent (see
\cite[Proposition 4.7]{QUEF-LI-LE-RO-PI}):
\begin{stprop}
\label{prop|nabla_0_implies_nabla_2}Let $\psi$ be an Orlicz function.
\begin{enumerate}
\item If $\psi$ satisfies the uniform $\nabla_{0}$-Condition, then it
satisfies the $\nabla_{2}$-Condition;
\item If $\psi$ satisfies the $\Delta^{2}$-Condition, then it satisfies
the uniform $\nabla_{0}$-Condition.
\end{enumerate}
\end{stprop}
\smallskip{}
For any $1<p<\infty$, every function $x\longmapsto x^{p}$ is an
Orlicz function which satisfies the uniform $\nabla_{0}$-Condition,
(so $\nabla_{2}$ and $\nabla_{0}$-conditions too) and the $\Delta_{2}$-Condition.
At the opposite side, for any $a>0$ and $b\geq1$, $x\longmapsto e^{ax^{b}}-1$
belongs to the $\Delta^{2}$-Class (and then to the uniform $\nabla_{0}$-Class),
yet not to the $\Delta_{2}$-one. In addition, the Orlicz functions
which can be written $x\rightarrow\exp\left(a\left(\ln\left(x+1\right)\right)^{b}\right)-1$
for $a>0$ and $b\geq1$, satisfy the $\nabla_{2}$ and $\nabla_{0}$-Conditions,
but do not belong to the $\Delta^{2}$-Class.

\smallskip{}
For a complete study of Orlicz spaces, we refer to \cite{KRASNO-RUT}
and \cite{RAO-REN}. We can also find precise and useful information
in \cite{QUEF-LI-LE-RO-PI}, such as other classes of Orlicz functions
and their links with each other.

\subsection{Background results}

All the results of the present paper are based on characterizations
of the boundedness and compactness of composition operators on Hardy-Orlicz
and Bergman-Orlicz spaces (\cite{Charp1,Charp2}). As I already said,
these characterizations essentially depend on the manner in which
the Orlicz function grows.

\smallskip{}

The characterizations of the boundedness and compactness of $C_{\phi}$
involve adapted Carleson measures, and then geometric notions. For
$\zeta\in\mathbb{S}_{N}$ and $0<h<1$, let us denote by $S\left(\zeta,h\right)$
and $\mathcal{S}\left(\zeta,h\right)$ the non-isotropic {}``balls'',
respectively in $\mathbb{B}_{N}$ and $\overline{\mathbb{B}_{N}}$,
defined by\[
S\left(\zeta,h\right)=\left\{ z\in\mathbb{B}_{N},\,\left|1-\left\langle z,\zeta\right\rangle \right|<h\right\} \mbox{ and }\mathcal{S}\left(\zeta,h\right)=\left\{ z\in\overline{\mathbb{B}_{N}},\,\left|1-\left\langle z,\zeta\right\rangle \right|<h\right\} .\]
\smallskip{}
We say that a finite positive Borel measure $\mu$ on $\overline{\mathbb{B}_{N}}$
is a $\psi$-Carleson measure, $\psi$ an Orlicz function, if\[
\mu\left(\mathcal{S}\left(\zeta,h\right)\right)=O_{h\rightarrow0}\left(\frac{1}{\psi\left(A\psi^{-1}\left(1/h^{N}\right)\right)}\right),\]
uniformly in $\zeta\in\mathbb{S}_{N}$ and for some constant $A>0$.
$\mu$ is a \emph{vanishing} $\psi$-Carleson measure if the above
condition is satisfied for every $A>0$ and with the big-Oh condition
replaced by a little-oh condition.

\smallskip{}
A finite positive Borel measure $\mu$ on $\mathbb{B}_{N}$ is a $\left(\psi,\alpha\right)$-Bergman-Carleson
measure if\[
\mu\left(S\left(\zeta,h\right)\right)=O_{h\rightarrow0}\left(\frac{1}{\psi\left(A\psi^{-1}\left(1/h^{N\left(\alpha\right)}\right)\right)}\right),\]
uniformly in $\zeta\in\mathbb{S}_{N}$ and for some constant $A>0$.
$\mu$ is a \emph{vanishing} $\left(\psi,\alpha\right)$-Bergman-Carleson
measure if the above condition is satisfied for every $A>0$ and with
the big-Oh condition replaced by a little-oh condition.

\smallskip{}

When $\psi$ satisfies the $\Delta_{2}$-Condition, a (vanishing)
$\psi$-Carleson measure (resp. (vanishing) $\left(\psi,\alpha\right)$-Bergman-Carleson
measure) is a (vanishing) Carleson measure (resp. (vanishing) Bergman-Carleson
measure) (see \cite[Sections 3]{Charp1,Charp2}).\smallskip{}

For $\phi:\mathbb{B}_{N}\rightarrow\mathbb{B}_{N}$, we denote by
$\mu_{\phi}$ the pull-back measure of $\sigma_{N}$ by the boundary
limit $\phi^{*}$ of $\phi$, and by $\mu_{\phi,\alpha}$ that of
$dv_{\alpha}$ by $\phi$. To be precise, for any $E\subset\overline{\mathbb{B}_{N}}$
(resp. $E\subset\mathbb{B}_{N}$),\[
\mu_{\phi}\left(E\right)=\sigma_{N}\left(\left(\phi^{*}\right)^{-1}\left(E\right)\right)\mbox{ and }\mu_{\phi,\alpha}\left(E\right)=v_{\alpha}\left(\phi^{-1}\left(E\right)\right).\]

\subsubsection{Results for Hardy-Orlicz spaces (see \cite[Section 3]{Charp2})}

The main theorem is the following:
\begin{stthm}
\label{thm|main_thm_H_O_rappel}Let $\psi$ be an Orlicz function
which satisfies the $\nabla_{2}$-Condition and let $\phi:\mathbb{B}_{N}\rightarrow\mathbb{B}_{N}$
be holomorphic.
\begin{enumerate}
\item If $\psi$ satisfies the uniform $\nabla_{0}$-Condition, then $C_{\phi}$
is bounded from $H^{\psi}\left(\mathbb{B}_{N}\right)$ into itself
if and only if $\mu_{\phi}$ is a $\psi$-Carleson measure.
\item If $\psi$ satisfies the $\nabla_{0}$-Condition, then $C_{\phi}$
is compact from $H^{\psi}\left(\mathbb{B}_{N}\right)$ into itself
if and only if $\mu_{\phi}$ is a vanishing $\psi$-Carleson measure.
\item If $\psi$ satisfies the $\Delta_{2}$-Condition, then $C_{\phi}$
is bounded (resp. compact) from $H^{\psi}\left(\mathbb{B}_{N}\right)$
into itself if and only if $\mu_{\phi}$ is a Carleson measure (resp.
a vanishing Carleson measure).
\item If $\psi$ satisfies the $\Delta^{2}$-Condition, then $C_{\phi}$
is bounded on $H^{\psi}\left(\mathbb{B}_{N}\right)$.
\end{enumerate}
\end{stthm}
The first two points are contained in \cite[Theorem 3.2]{Charp2};
according to \cite[Theorem 3.35]{COWEN-MACCLUER}, the third point
means that, if $\psi$ satisfies the $\Delta_{2}$-Condition, then
$C_{\phi}$ is bounded (resp. compact) on $H^{\psi}\left(\mathbb{B}_{N}\right)$
if and only if it is on $H^{p}\left(\mathbb{B}_{N}\right)$ (see \cite[Corollary 3.4]{Charp2}).
The last point is \cite[Theorem 3.7]{Charp2}.\smallskip{}

Due to the non-separability of small Hardy-Orlicz spaces, \cite[Theorem 3.2]{Charp2}
is not a direct consequence of Carleson-type embedding theorems obtained
in \cite[Section 2]{Charp2}; however, if we follow the proofs of
these embedding theorems directly for composition operators, by working
on spheres of radius $0<r<1$, then we get the following characterizations
of both boundedness and compactness of composition operators:
\begin{stthm}
\label{thm|cond_nec_comp_C_phi_H_O}Let $\psi$ be an Orlicz function
satisfying the $\nabla_{2}$-Condition and let $\phi:\mathbb{B}_{N}\rightarrow\mathbb{B}_{N}$
be holomorphic.
\begin{enumerate}
\item If $\psi$ satisfies the uniform $\nabla_{0}$-Condition, then $C_{\phi}$
is bounded on $H^{\psi}\left(\mathbb{B}_{N}\right)$ if and only if
there exists some $A>0$ such that\begin{equation}
\sup_{0<r<1}\mu_{\phi_{r}}\left(\mathcal{S}\left(\zeta,h\right)\right)=O_{h\rightarrow0}\left(\frac{1}{\psi\left(A\psi^{-1}\left(1/h^{N}\right)\right)}\right)\label{eq|carac_bis_cont_H_O}\end{equation}
uniformly in $\zeta\in\mathbb{S}_{N}$.
\item If $\psi$ satisfies the $\nabla_{0}$-Condition, then $C_{\phi}$
is compact on $H^{\psi}\left(\mathbb{B}_{N}\right)$ if and only if,
for every $A>0$,\begin{equation}
\sup_{0<r<1}\mu_{\phi_{r}}\left(\mathcal{S}\left(\zeta,h\right)\right)=O_{h\rightarrow0}\left(\frac{1}{\psi\left(A\psi^{-1}\left(1/h^{N}\right)\right)}\right)\label{eq|carac_bis_compa_H_O}\end{equation}
uniformly in $\zeta\in\mathbb{S}_{N}$.
\end{enumerate}
\end{stthm}
In the further, we will see how these two characterizations are useful
depending on the situations.

\subsubsection{Results for Bergman-Orlicz spaces (see \cite[Section 3]{Charp1})}

The main result is similar to that stated in the previous paragraph:

\begin{stthm}
\label{thm|main_thm_B_O_rappel}Let $\psi$ be an Orlicz function,
let $\alpha>-1$ and let $\phi:\mathbb{B}_{N}\rightarrow\mathbb{B}_{N}$
be holomorphic.
\begin{enumerate}
\item If $\psi$ satisfies the uniform $\nabla_{0}$-Condition, then $C_{\phi}$
is bounded from $A_{\alpha}^{\psi}\left(\mathbb{B}_{N}\right)$ into
itself if and only if $\mu_{\phi}$ is a $\left(\psi,\alpha\right)$-Bergman-Carleson
measure.
\item If $\psi$ satisfies the $\nabla_{0}$-Condition, then $C_{\phi}$
is compact from $A_{\alpha}^{\psi}\left(\mathbb{B}_{N}\right)$ into
itself if and only if $\mu_{\phi}$ is a vanishing $\left(\psi,\alpha\right)$-Bergman-Carleson
measure.
\item If $\psi$ satisfies the $\Delta_{2}$-Condition, then $C_{\phi}$
is bounded (resp. compact) from $A_{\alpha}^{\psi}\left(\mathbb{B}_{N}\right)$
into itself if and only if $\mu_{\phi}$ is a Bergman-Carleson measure
(resp. a vanishing Bergman-Carleson measure).
\item If $\psi$ satisfies the $\Delta^{2}$-Condition, then $C_{\phi}$
is bounded on $A_{\alpha}^{\psi}\left(\mathbb{B}_{N}\right)$.
\end{enumerate}
\end{stthm}
According to \cite[Theorem 3.37]{COWEN-MACCLUER}, the third point
means that, if $\psi$ satisfies the $\Delta_{2}$-Condition, then
$C_{\phi}$ is bounded (resp. compact) on $A_{\alpha}^{\psi}\left(\mathbb{B}_{N}\right)$,
if and only if $C_{\phi}$ is bounded (resp. compact) on $A_{\alpha}^{p}\left(\mathbb{B}_{N}\right)$.

\section{Main results}

\subsection{Compactness of $C_{\phi}$ on every Hardy-Orlicz or Bergman-Orlicz
spaces}

The following theorem completes both \prettyref{thm|main_thm_H_O_rappel}
and \prettyref{thm|main_thm_B_O_rappel}.
\begin{stthm}
\label{thm|theorem_Comp_on_every_H_O_H_infinity}Let $\phi:\mathbb{B}_{N}\rightarrow\mathbb{B}_{N}$
be a holomorphic map. The following assertions are equivalent:
\begin{enumerate}
\item $C_{\phi}$ is compact on $H^{\psi}\left(\mathbb{B}_{N}\right)$,
for every Orlicz function $\psi$;
\item For some $\alpha>-1$, $C_{\phi}$ is compact on $A_{\alpha}^{\psi}\left(\mathbb{B}_{N}\right)$,
for every Orlicz function $\psi$;
\item $C_{\phi}$ is compact on $H^{\infty}\left(\mathbb{B}_{N}\right)$;
\item $\left\Vert \phi\right\Vert _{\infty}<1$.
\end{enumerate}
\end{stthm}
\begin{proof}
It is well-known that $C_{\phi}$ is compact on $H^{\infty}\left(\mathbb{B}_{N}\right)$
if and only if $\left\Vert \phi\right\Vert _{\infty}<1$. Using the
fact that a composition operator is compact on $H^{\psi}\left(\mathbb{B}_{N}\right)$
(resp. $A_{\alpha}^{\psi}\left(\mathbb{B}_{N}\right)$) if and only
if for every bounded sequence $\left(f_{n}\right)_{n}\subset H^{\psi}\left(\mathbb{B}_{N}\right)$,
$\left\Vert f_{n}\right\Vert _{\psi}\leq1$, (resp. $\left(f_{n}\right)_{n}\subset A_{\alpha}^{\psi}\left(\mathbb{B}_{N}\right)$,
$\left\Vert f_{n}\right\Vert _{A_{\alpha}^{\psi}}\leq1$) which tends
to $0$ uniformly on every compact subset of $\mathbb{B}_{N}$, then
$\left\Vert f_{n}\circ\phi\right\Vert _{\psi}\xrightarrow[n\rightarrow\infty]{}0$
(resp. $\left\Vert f_{n}\circ\phi\right\Vert _{A_{\alpha}^{\psi}}\xrightarrow[n\rightarrow\infty]{}0$),
it is not difficult to show that if $\left\Vert \phi\right\Vert _{\infty}<1$,
then $C_{\phi}$ is compact on $H^{\psi}\left(\mathbb{B}_{N}\right)$
(resp. $A_{\alpha}^{\psi}\left(\mathbb{B}_{N}\right)$.)

It remains to prove (1)$\Rightarrow$(4) and (2)$\Rightarrow$(4).
We first deal with the proof of (1)$\Rightarrow$(4). We will use
the necessary part of \prettyref{thm|cond_nec_comp_C_phi_H_O}. Let
us assume that $\phi$ induces a compact composition operator on every
Hardy-Orlicz space. According to \prettyref{eq|carac_bis_compa_H_O},
this means that\[
\sup_{0<r<1}\sup_{\zeta\in\mathbb{S}_{N}}\left(\mu_{\phi_{r}}\left(\mathcal{S}\left(\zeta,h\right)\right)\right)=o_{h\rightarrow0}\left(\frac{1}{\psi\left(A\psi^{-1}\left(1/h^{N}\right)\right)}\right),\]
for every $A>0$ and every Orlicz function $\psi$, which in turn
implies\begin{equation}
\sup_{0<r<1}\sup_{\zeta\in\mathbb{S}_{N}}\left(\mu_{\phi_{r}}\left(\mathcal{S}\left(\zeta,h\right)\right)\right)\leq\frac{1}{\psi\left(A\psi^{-1}\left(1/h^{N}\right)\right)},\label{eq|eq_1_comp_on_every_H_O}\end{equation}
for every $A>0$, for every Orlicz function $\psi$ and for $h$ sufficiently
small. We intend to show that \[
\sup_{0<r<1}\sup_{\zeta\in\mathbb{S}_{N}}\left(\mu_{\phi_{r}}\left(\mathcal{S}\left(\zeta,h\right)\right)\right)=0,\]
for all $0<h\leq h_{0}$, $h_{0}\in\left(0,1\right)$. By contradiction,
we suppose that ${\displaystyle \sup_{0<r<1}\sup_{\zeta\in\mathbb{S}_{N}}\left(\mu_{\phi_{r}}\left(\mathcal{S}\left(\zeta,h\right)\right)\right)\neq0}$
for every $h>0$, since \[
h\longmapsto\sup_{0<r<1}\sup_{\zeta\in\mathbb{S}_{N}}\left(\mu_{\phi_{r}}\left(\mathcal{S}\left(\zeta,h\right)\right)\right)\]
is an increasing function on $\left(0,1\right)$. A straightforward
computation shows that inequality \prettyref{eq|eq_1_comp_on_every_H_O}
is satisfied for every $A>0$, for every Orlicz function $\psi$ and
for $h$ small enough, if and only if we have, by putting $x=1/h$,\begin{equation}
\frac{\psi^{-1}\left(x^{N}\right)}{\psi^{-1}\left({\displaystyle 1/{\displaystyle \sup_{0<r<1}\sup_{\zeta\in\mathbb{S}_{N}}\left(\mu_{\phi_{r}}\left(\mathcal{S}\left(\zeta,1/x\right)\right)\right)}}\right)}\leq\frac{1}{A},\label{eq|eq_2_comp_on_every_H_O}\end{equation}
for every $A>0$, for every Orlicz function $\psi$ and for $x$ large
enough. The following lemma ensures that this cannot occur:
\begin{stlem}
\label{lem|concave}Let $f,g:\left[0,+\infty\right[\rightarrow\left[0,+\infty\right[$
be two increasing functions which tend to $+\infty$ at $+\infty$.
There exist $\delta>0$ and a continuous increasing concave function
$\nu:\left[0,+\infty\right[\rightarrow\left[0,+\infty\right[$, with
${\displaystyle \lim_{x\rightarrow+\infty}\nu\left(x\right)=+\infty}$,
such that ${\displaystyle \frac{\nu\left(f\left(x\right)\right)}{\nu\left(g\left(x\right)\right)}\geq\delta>0}$,
for every $x$ large enough.
\end{stlem}
We assume for a while that this lemma has been proven, and we finish
the proof of \prettyref{thm|theorem_Comp_on_every_H_O_H_infinity}.
With the notations of the lemma, we put\[
\left\{ \begin{array}{l}
f\left(x\right)=x^{N}\mbox{ and}\\
g\left(x\right){\textstyle =1/{\displaystyle \sup_{0<r<1}\sup_{\zeta\in\mathbb{S}_{N}}\left(\mu_{\phi_{r}}\left(\mathcal{S}\left(\zeta,1/x\right)\right)\right)}}.\end{array}\right.\]
It is clear that $\lim_{x\rightarrow+\infty}f\left(x\right)=+\infty$;
since $C_{\phi}$ is supposed to be compact on every $H^{\psi}\left(\mathbb{B}_{N}\right)$,
it is in particular bounded on $H^{p}\left(\mathbb{B}_{N}\right)$
(\cite[Corollary 3.5]{Charp2}), then we have $g\left(x\right)\xrightarrow[x\rightarrow+\infty]{}+\infty$
(\prettyref{thm|main_thm_H_O_rappel}, 3). Now, the above lemma provides
a constant $\delta>0$ and a continuous increasing concave function
$\nu$, tending to infinity at infinity, such that\[
\nu\left(x^{N}\right)/\nu\left({\displaystyle \frac{1}{{\displaystyle \sup_{\xi\in\mathbb{S}_{N}}\left(\mu_{\phi}\left(S_{f}\left(\xi,1/x\right)\right)\right)}}}\right)\geq\delta>0,\]
for every $x$ large enough. It is not difficult to check that $\nu$
can be constructed such that $\psi=\nu^{-1}$ is an Orlicz function,
i.e. such that ${\displaystyle \frac{x}{\nu\left(x\right)}\xrightarrow[x\rightarrow+\infty]{}+\infty}$
(that is what we are doing in the proof of the lemma below). Therefore,
we get a contradiction with Condition \prettyref{eq|eq_2_comp_on_every_H_O},
so we must have\[
\sup_{0<r<1}\sup_{\xi\in\mathbb{S}_{N}}\left(\mu_{\phi_{r}}\left(S_{f}\left(\xi,h\right)\right)\right)=0,\]
for every $h>0$ small enough. It follows that there exists some $0<r_{0}<1$
such that\begin{equation}
\sup_{0<r<1}\mu_{\phi_{r}}\left(\mathcal{C}\left(r_{0},1\right)\right)=0,\label{eq|eq_3_comp_on_every_H_psi}\end{equation}
where $\mathcal{C}\left(r_{0},1\right)=\left\{ z\in\mathbb{B}_{N},\, r_{0}<\left|z\right|<1\right\} $.
We intend to show that $\phi^{-1}\left(\mathcal{C}\left(r_{0},1\right)\right)=\emptyset$,
which should give the result. Let $0<r<1$ and let us look at the
set\[
\phi_{r}^{-1}\left(\mathcal{C}\left(r_{0},1\right)\right)\cap\mathbb{S}_{N}=\left\{ \zeta\in\mathbb{S}_{N},\,\phi_{r}\left(\zeta\right)\in\mathcal{C}\left(r_{0},1\right)\right\} .\]
Condition \prettyref{eq|eq_3_comp_on_every_H_psi} implies\[
\sigma_{N}\left(\phi_{r}^{-1}\left(\mathcal{C}\left(r_{0},1\right)\right)\cap\mathbb{S}_{N}\right)=0.\]
Since $\phi_{r}$ is continuous on $\overline{\mathbb{B}_{N}}$, $\phi_{r}^{-1}\left(\mathcal{C}\left(r_{0},1\right)\right)\cap\mathbb{S}_{N}$
must be an open subset of $\mathbb{S}_{N}$ and then must be empty.
So we have proven that, for any $r\in\left(0,1\right)$,\[
\left\{ \zeta\in\mathbb{S}_{N},\,\phi_{r}\left(\zeta\right)\in\mathcal{C}\left(r_{0},1\right)\right\} =\phi^{-1}\left(\mathcal{C}\left(r_{0},1\right)\right)\cap r\mathbb{S}_{N}=\emptyset,\]
where $r\mathbb{S}_{N}=\left\{ z\in\overline{\mathbb{B}_{N}},\,\left|z\right|=r\right\} $,
hence \[
\phi^{-1}\left(\mathcal{C}\left(r_{0},1\right)\right)=\bigcup_{0<r<1}\left(\phi^{-1}\left(\mathcal{C}\left(r_{0},1\right)\right)\cap r\mathbb{S}_{N}\right)=\emptyset.\]

The proof in the Bergman-Orlicz context is much easier. Proceeding
as above and using the necessary part of the second point of \prettyref{thm|main_thm_B_O_rappel},
we get that Condition $\mu_{\phi}\left(\mathcal{C}\left(r_{0},1\right)\right)=0$
must hold, for some $0<r_{0}<1$. By continuity of the map $\phi$
on $\mathbb{B}_{N}$, $\phi^{-1}\left(\mathcal{C}\left(r_{0},1\right)\right)$
cannot be but empty.

\smallskip{}

To finish the proof, we have to do that of \prettyref{lem|concave}:
\begin{proof}
[Proof of \prettyref{lem|concave}]The proof will be constructive.
Let $f$ and $g$ be given as in the statement of the lemma. We are
going to build by induction a sequence $\left(a_{n}\right)_{n}$ which
will be of interest in the construction of the desired function $\nu$.
We put $a_{0}=0$, $a_{1}=1$, and we deduce $a_{n+2}$ from $a_{n}$
and $a_{n+1}$ in the following way: we define\[
b_{n+2}=\sup\left\{ g\left(x\right),\, f\left(x\right)\leq a_{n+1}\right\} \]
and\[
a_{n+2}=\max\left\{ b_{n+2},a_{n+1}+\left(a_{n+1}-a_{n}\right)\right\} .\]
We observe that:
\begin{enumerate}
\item If $f\left(x\right)\leq a_{n+1}$, then $g\left(x\right)\leq a_{n+2}$;
\item $a_{n+2}-a_{n+1}\geq a_{n+1}-a_{n}\geq1.$
\end{enumerate}
We now construct the concave function $\nu$ as a continuous affine
one, whose derivative is equal to ${\displaystyle \varepsilon_{n}=\frac{1}{\sqrt{n}\left(a_{n+1}-a_{n}\right)}}$
on the interval $\left(a_{n},a_{n+1}\right)$, and with $\nu\left(0\right)=0$.
Of course $\nu$ is increasing and then maps $\left[0,+\infty\right[$
into itself. Since $\varepsilon_{n}$ is decreasing, because of 2)
above, $\nu$ is concave. In order to check that $\nu$ tends to infinity
at infinity, we compute $\nu\left(a_{n}\right)$:\[
\nu\left(a_{n+1}\right)=\nu\left(a_{n}\right)+\varepsilon_{n}\left(a_{n+1}-a_{n}\right)=\nu\left(a_{n}\right)+\frac{1}{\sqrt{n}}.\]
Therefore $\nu\left(a_{n+1}\right)=\sum_{k=1}^{n+1}\frac{1}{\sqrt{k}}$
which shows that $\lim_{x\rightarrow+\infty}\nu\left(x\right)=+\infty$,
since $a_{n}\rightarrow+\infty$.

We now check that ${\displaystyle \frac{\nu\circ f\left(x\right)}{\nu\circ g\left(x\right)}}$
is bounded below by some constant $\delta>0$, when $x$ is big enough.
Let $x\in\left[0,+\infty\right[$, and let $n$ be an integer such
that $a_{n}\leq f\left(x\right)\leq a_{n+1}$; we have $\nu\left(f\left(x\right)\right)\geq\nu\left(a_{n}\right)$.
Using the first property of the sequence $\left(a_{n}\right)_{n}$
above, we get $\nu\left(g\left(x\right)\right)\leq\nu\left(a_{n+2}\right)$.
This yields, for $n\geq1$,\[
\frac{\nu\left(f\left(x\right)\right)}{\nu\left(g\left(x\right)\right)}\geq\frac{\nu\left(a_{n}\right)}{\nu\left(a_{n+2}\right)}=\frac{\sum_{k=1}^{n}\frac{1}{\sqrt{k}}}{\sum_{k=1}^{n+2}\frac{1}{\sqrt{k}}}\geq\delta>0,\]
hence the result.
\end{proof}
\end{proof}

\subsection{Kor\'anyi regions and compactness of $C_{\phi}$ on Hardy-Orlicz
and Bergman-Orlicz spaces}

For $\zeta\in\mathbb{S}_{N}$ and $a>1$, we recall that the Kor\'anyi
approach region $\Gamma\left(\zeta,a\right)$ of angular opening $a$
is defined by\[
\Gamma\left(\zeta,a\right)=\left\{ z\in\mathbb{B}_{N},\,\left|1-\left\langle z,\zeta\right\rangle \right|<\frac{a}{2}\left(1-\left|z\right|^{2}\right)\right\} .\]

\cite[Theorem 6.4]{COWEN-MACCLUER} and the third part of \prettyref{thm|main_thm_H_O_rappel}
yields the following result:
\begin{stthm}
\label{thm|thm_Kor_positive_Delta_2}Let $\psi$ be an Orlicz function
satisfying the $\Delta_{2}\cap\nabla_{2}$-Condition. Let also $\phi:\mathbb{B}_{N}\rightarrow\mathbb{B}_{N}$
be holomorphic. We assume that $N>1$ and we fix $b_{N}=\left(\cos\left(\pi/\left(2N\right)\right)\right)^{-1}$.
\begin{enumerate}
\item If $\phi\left(\mathbb{B}_{N}\right)\subset\Gamma\left(\zeta,b_{N}\right)$,
then $C_{\phi}$ is bounded on $H^{\psi}\left(\mathbb{B}_{N}\right)$;
\item If $\phi\left(\mathbb{B}_{N}\right)\subset\Gamma\left(\zeta,b\right)$,
for some $\zeta\in\mathbb{S}_{N}$ and for some $1<b<b_{N}$, then
$C_{\phi}$ is compact on $H^{\psi}\left(\mathbb{B}_{N}\right)$.
\item Both of the above results are sharp in the following sense: given
$c>b_{N}$, there exists $\phi$ with $\phi\left(\mathbb{B}_{N}\right)\subset\Gamma\left(\zeta,c\right)$,
for some $\zeta\in\mathbb{S}_{N}$, and $C_{\phi}$ not bounded on
$H^{\psi}\left(\mathbb{B}_{N}\right)$; there also exists some $\phi$
with $\phi\left(\mathbb{B}_{N}\right)\subset\Gamma\left(\zeta,b_{N}\right)$,
for some $\zeta\in\mathbb{S}_{N}$, and $C_{\phi}$ not compact on
$H^{\psi}\left(\mathbb{B}_{N}\right)$.
\end{enumerate}
\end{stthm}
\begin{strem}
(1) If $N=1$, the two first point of the previous theorem are true
if we put $b_{1}=+\infty$ and $\Gamma\left(\zeta,+\infty\right)=\mathbb{D}$.
Indeed, the first point is nothing but the continuity of every composition
operator on the disc, and the second one is contained in \cite[Proposition 3.25]{COWEN-MACCLUER}
which says that, whenever $\phi\left(\mathbb{D}\right)$ is contained
in some nontangential approach region in $\mathbb{D}$, then $C_{\phi}$
is Hilbert-Schmidt in $H^{2}\left(\mathbb{D}\right)$, hence compact
on every $H^{p}\left(\mathbb{D}\right)$, $1\leq p<\infty$.

(2) Following the proof of \cite[Theorem 3.3]{CLAHANE}, it is not
difficult to show that the boundedness or compactness of $C_{\phi}$
on $H^{\psi}\left(\mathbb{B}_{N}\right)$ implies that on $A_{\alpha}^{\psi}\left(\mathbb{B}_{N}\right)$,
for any $\alpha>-1$, as soon as the Orlicz function $\psi$ satisfies
the $\Delta_{2}$-Condition. Thus, the first two points of the previous
theorem also holds for Bergman-Orlicz spaces.
\end{strem}
The following result shows that \prettyref{thm|thm_Kor_positive_Delta_2}
does not hold as soon as the Orlicz function grows fast.
\begin{stthm}
\label{thm|thm_pas_comp_delta_2_Kor_app}Let $\psi$ be an Orlicz
function satisfying the $\Delta^{2}$-Condition. Then, for every $\zeta\in\mathbb{S}_{N}$
and every $b>1$, there exists a holomorphic self-map $\phi$ taking
$\mathbb{B}_{N}$ into $\Gamma\left(\zeta,b\right)$, such that $C_{\phi}$
is not compact on $H^{\psi}\left(\mathbb{B}_{N}\right)$.\end{stthm}
\begin{strem}
Observe that there is no assumption on $N$.\end{strem}
\begin{proof}
[Proof of \prettyref{thm|thm_pas_comp_delta_2_Kor_app}]The proof
will use the necessary part of the second point of \prettyref{thm|main_thm_H_O_rappel}.
First of all, we recall that $\Delta^{2}$-Condition implies $\nabla_{2}$-Condition
(see Paragraph 2.2). We denote by $e_{1}$ the vector $\left(1,0\ldots0\right)$
in $\mathbb{C}^{N}$. It is clearly enough to prove the theorem for
$\zeta=e_{1}$. For any $b>1$, we set\begin{equation}
\beta=\frac{2\cos^{-1}\left(1/b\right)}{\pi}\label{eq|def_beta}\end{equation}
in $\left(0,1\right)$. We need a lemma whose proof is included in
that of \cite[Theorem 6.4]{COWEN-MACCLUER}:
\begin{stlem}
\label{lem|lemma_1_prop_Kor}Let $b>1$ and let $\beta$ be defined
by \prettyref{eq|def_beta}. There exists a holomorphic map $\phi:\mathbb{B}_{N}\rightarrow\mathbb{B}_{N}$,
with $\phi\left(\mathbb{B}_{N}\right)\subset\Gamma\left(e_{1},b\right)$,
such that\begin{equation}
\sigma_{N}\phi^{-1}\left(\mathcal{S}\left(e_{1},h\right)\right)\geq Ch^{1/\beta},\label{eq|lemme1_lower_estimate}\end{equation}
for some constant $C>0$ depending only on $\phi$ and $b$.\end{stlem}
\begin{proof}
Without going into details, we briefly give the ideas of the proof.
It uses the deep Alexandrov's result which gives the existence of
non-constant inner functions in $\mathbb{B}_{N}$ (\cite{ALEXANDROV}).
Therefore, we consider a function $\phi$ which can be written\[
\phi=\left(\kappa\circ\varphi,0'\right),\]
where $0'$ is the $\left(n-1\right)$-tuple $\left(0,\ldots,0\right)$,
$\varphi$ is an inner function with $\varphi\left(0\right)=0$, and
where $\kappa$ is a biholomorphic map from $\mathbb{D}$ onto the
non-tangential approach region $\Gamma\left(1,b\right)$ in the disc,
defined by\[
\Gamma\left(1,b\right)=\left\{ z\in\mathbb{D},\left|1-z\right|<\frac{b}{2}\left(1-\left|z\right|^{2}\right)\right\} .\]
One can show that the lower-estimate \prettyref{eq|lemme1_bis_lower_estimate}
holds for this map $\phi$, using the fact that inner functions $\varphi$
are measure preserving maps of $\mathbb{S}_{N}$ into $\mathbb{T}$
(see \cite[p. 405]{RUDIN_80}) in the following sense:\[
\sigma_{N}\left(\left(\varphi^{*}\right)^{-1}\left(E\right)\right)=\sigma_{1}\left(E\right),\]
for any Borel set $E$ in $\mathbb{T}$.
\end{proof}
Let $\phi$ be as in the statement of the theorem. According to the
necessary part of the second point of \prettyref{thm|main_thm_H_O_rappel},
the previous lemma ensures that, if we show that for any Orlicz function
$\psi$ satisfying the $\Delta^{2}$-Condition, for any $\beta\in\left(0,1\right)$,
there exists some $A>0$ such that\begin{equation}
\frac{1}{\psi\left(A\psi^{-1}\left(1/h^{N}\right)\right)}\leq h^{1/\beta},\label{eq|eq_thm_MAIN_2_proof}\end{equation}
for every $h$ small enough, then $C_{\phi}$ would not be compact
on $H^{\psi}\left(\mathbb{B}_{N}\right)$. Now, putting $y=\psi^{-1}\left(1/h^{N}\right)$,
an easy computation implies that \prettyref{eq|eq_thm_MAIN_2_proof}
is equivalent to ${\displaystyle \psi\left(y\right)^{\frac{1}{N\beta}}}\leq\psi\left(Cy\right)$
for some constant $C>0$. We conclude the proof by noticing that this
latter condition is trivial if ${\displaystyle 0<\frac{1}{N\beta}\leq1}$
(which allows to recover the compactness part of the third point of
\prettyref{thm|thm_Kor_positive_Delta_2}), while it is nothing but
$\Delta^{2}$-Condition if ${\displaystyle \frac{1}{N\beta}>1}$ (see
Paragraph 2.2).\end{proof}
\begin{strem}
When $N=1$, the proof of \prettyref{lem|lemma_1_prop_Kor} can be
simplified: first, because the existence of a non-constant inner function
in the unit disc is trivial, and then because it clearly suffices
to take $\varphi\left(z\right)=z$, what just turns the proof of \prettyref{lem|lemma_1_prop_Kor}
into considering a biholomorphic map $\kappa$ from $\mathbb{D}$
onto an non-tangential approach region.\smallskip{}

\end{strem}
The previous remark leads us to say some words about weighted Bergman-Orlicz
spaces in dimension one. Indeed, we can adapt the proof of \prettyref{lem|lemma_1_prop_Kor}
to get the following result.
\begin{stlem}
\label{lem|Lemme_1__Kor_bis_B_O}Let $\alpha>-1$, let $b>1$ and
let $\beta$ be defined by \prettyref{eq|def_beta}. There exists
a holomorphic map $\phi:\mathbb{D}\rightarrow\mathbb{D}$, $\phi\left(\mathbb{D}\right)\subset\Gamma\left(1,b\right)$,
such that\begin{equation}
v_{\alpha}\left(\phi^{-1}\left(S\left(1,h\right)\right)\right)\geq Ch^{\left(2+\alpha\right)/\beta},\label{eq|lemme1_bis_lower_estimate}\end{equation}
for some constant $C>0$ depending only on $\alpha$, $\phi$ and
$b$.
\end{stlem}
For the seek of completeness, we prefer to give some details of the
proof of this lemma, in order to point out the slightly difference
with that of \cite[Theorem 6.4, 3)]{COWEN-MACCLUER}.
\begin{proof}
We consider a biholomorphic map $\kappa$ from $\mathbb{D}$ onto\[
{\displaystyle \Gamma\left(1,b\right)=\left\{ z\in\mathbb{D},\,\left|1-z\right|<\frac{b}{2}\left(1-\left|z\right|^{2}\right)\right\} },\]
for some $b>1$. As it is explained in the proof of \cite[Theorem 6.4, 3)]{COWEN-MACCLUER},
if $\beta\in\left(0,1\right)$ is defined by \prettyref{eq|def_beta},
then the function ${\displaystyle g\left(z\right):=\frac{1-\kappa\left(z\right)}{\left(1-z\right)^{\beta}}}$
is continuous and non-zero in $\overline{\mathbb{D}}\cap V$, where
$V$ is a closed disc (with non-empty interior) centered at $1$,
and $\kappa^{-1}\left(S\left(1,h\right)\right)\subset\overline{\mathbb{D}}\cap V$,
for $h>0$ sufficiently small. Then, for such $h$, we follow the
computation at the end of the proof of \cite[Theorem 6.4, 3)]{COWEN-MACCLUER}
to get \[
\kappa^{-1}\left(S\left(1,h\right)\right)\supset S\left(1,\tilde{C}h^{1/\beta}\right),\]
for some constant $\tilde{C}>0$, depending only on $\kappa$ and
$b$. Therefore,\[
v_{\alpha}\left(\kappa^{-1}\left(S\left(1,h\right)\right)\right)\geq Ch^{\left(2+\alpha\right)/\beta},\]
where $C>0$ depends on $\alpha$, $\kappa$ and $b$.
\end{proof}
Now, it is sufficient to argue as in in the end of the proof of \prettyref{thm|thm_pas_comp_delta_2_Kor_app}
to get:
\begin{stprop}
Let $\alpha>-1$, let $b>1$ and let $\psi$ be an Orlicz function
satisfying the $\Delta^{2}$-Condition. There exists a holomorphic
map $\phi:\mathbb{D}\rightarrow\mathbb{D}$, with range contained
in some non-tangential approach region $\Gamma\left(\zeta,b\right)$,
$\zeta\in\mathbb{T}$, such that the induced composition operator
$C_{\phi}$ is not compact on $A_{\alpha}^{\psi}\left(\mathbb{D}\right)$.
\end{stprop}
\smallskip{}

The proof of the previous proposition does not work directly when
$N>1$, because we do not know if there exists a non-constant inner
function which is measure-preserving from $\mathbb{B}_{N}$ to $\mathbb{D}$
in the following sense:\[
v_{\alpha}\left(\phi^{-1}\left(E\right)\right)=A_{\alpha}\left(E\right),\]
for any $E\subset\mathbb{D}$, where $A_{\alpha}$ is the weighted
area measure in $\mathbb{D}$.

\subsection{Another characterization of the compactness of $C_{\phi}$ on weighted
Bergman-Orlicz spaces}

The following result generalizes that obtained in \cite{ZHU_COMP_OP_BERG_BALL}
for classical Bergman spaces:
\begin{stthm}
\label{thm|thm_MAIN_3_B_O}We assume that $\alpha>-1$. Let $\phi:\mathbb{B}_{N}\rightarrow\mathbb{B}_{N}$
be holomorphic and let $\psi$ be an Orlicz function which satisfies
the $\nabla_{0}$-Condition. We assume that $C_{\phi}$ is bounded
from $A_{\beta}^{\psi}\left(\mathbb{B}_{N}\right)$ into itself for
some $-1<\beta<\alpha$. Then $C_{\phi}$ is compact from $A_{\alpha}^{\psi}\left(\mathbb{B}_{N}\right)$
into itself if and only if\begin{equation}
\lim_{\left|z\right|\rightarrow1}\frac{\psi^{-1}\left(1/\left(1-\left|\phi\left(z\right)\right|\right)^{N\left(\alpha\right)}\right)}{\psi^{-1}\left(1/\left(1-\left|z\right|\right)^{N\left(\alpha\right)}\right)}=0.\label{eq|cond_main_thm}\end{equation}
\end{stthm}
\begin{proof}
The proof of the necessary part is the same as that of \cite[Theorem 5.7]{QUEF-LI-LE-RO-PI}.
We deal with the proof of the sufficiency of \prettyref{eq|cond_main_thm}.
Without loss of generality, we assume that $\phi\left(0\right)=0$.
According to the second point of \prettyref{thm|main_thm_B_O_rappel},
by the convexity of the Orlicz function $\psi$, it is sufficient
to show that for every $B>0$, there exists $h_{0}\in\left(0,1\right)$,
such that\begin{equation}
\psi^{-1}\left(1/\mu_{\phi,\alpha}\left(S\left(\xi,h\right)\right)\right)\geq B\psi^{-1}\left(1/h^{N\left(\alpha\right)}\right),\label{eq|eq_0_main_thm_B_O}\end{equation}
uniformly in $\xi\in\mathbb{S}_{N}$, and for any $0<h<h_{0}$. Let
$\alpha$ and $\beta$ be as in the statement of the theorem. We have\begin{eqnarray}
\mu_{\phi,\alpha}\left(S\left(\xi,h\right)\right) & = & \int_{\phi^{-1}\left(S\left(\xi,h\right)\right)}\left(1-\left|z\right|^{2}\right)^{\alpha}dv\left(z\right)\nonumber \\
 & \leq & 2^{\alpha-\beta}\sup_{z\in\phi^{-1}\left(S\left(\xi,h\right)\right)}\left(1-\left|z\right|\right)^{\alpha-\beta}\int_{\phi^{-1}\left(S\left(\xi,h\right)\right)}\left(1-\left|z\right|^{2}\right)^{\beta}dv\left(z\right)\nonumber \\
 & = & 2^{\alpha-\beta}\sup_{z\in\phi^{-1}\left(S\left(\xi,h\right)\right)}\left(1-\left|z\right|\right)^{\alpha-\beta}\mu_{\phi,\beta}\left(S\left(\xi,h\right)\right)\nonumber \\
 & \leq & 2^{\alpha-\beta}\sup_{z\in\phi^{-1}\left(S\left(\xi,h\right)\right)}\left(1-\left|z\right|\right)^{\alpha-\beta}\frac{1}{\psi\left(C_{\beta}\psi^{-1}\left(1/h^{N\left(\beta\right)}\right)\right)},\label{eq|eq_2_main_thm_B_O}\end{eqnarray}
where the last inequality stands for some constant $C_{\beta}\geq1$
and for $h$ small enough, since $C_{\phi}$ is supposed to be bounded
on $A_{\beta}^{\psi}\left(\mathbb{B}_{N}\right)$.

Now, since $\alpha-\beta>0$, the hypothesis \prettyref{eq|cond_main_thm}
is equivalent to the fact that, for any $A>0$,\[
\left(1-\left|z\right|\right)^{\alpha-\beta}\leq\frac{1}{\left(\psi\left(A\psi^{-1}\left(1/\left(1-\left|\phi\left(z\right)\right|\right)^{N\left(\alpha\right)}\right)\right)\right)^{\frac{\alpha-\beta}{N\left(\alpha\right)}}},\]
whenever $\left|z\right|$ is close enough to $1$. Moreover, let
us observe that if $z\in\phi^{-1}\left(S\left(\xi,h\right)\right)$,
then\[
1-\left|z\right|\leq1-\left|\phi\left(z\right)\right|\leq\left|1-\left\langle \phi\left(z\right),\xi\right\rangle \right|<h\]
so that, for any $A>0$,\[
\sup_{z\in\phi^{-1}\left(S\left(\xi,h\right)\right)}\left(1-\left|z\right|\right)^{\alpha-\beta}\leq\frac{1}{\left(\psi\left(A\psi^{-1}\left(1/h^{N\left(\alpha\right)}\right)\right)\right)^{\frac{\alpha-\beta}{N\left(\alpha\right)}}},\]
for any $h>0$ small enough, using the fact that $\psi$ is a non-decreasing
function and that $\alpha-\beta>0$. Thus, it follows from \prettyref{eq|eq_2_main_thm_B_O}
that\[
\mu_{\phi,\alpha}\left(S\left(\xi,h\right)\right)\leq2^{\alpha-\beta}\frac{1}{\left(\psi\left(A\psi^{-1}\left(1/h^{N\left(\alpha\right)}\right)\right)\right)^{\frac{\alpha-\beta}{N\left(\alpha\right)}}}\frac{1}{\psi\left(C_{\beta}\psi^{-1}\left(1/h^{N\left(\beta\right)}\right)\right)},\]
for any $A>0$ and $h$ small enough. Using \prettyref{eq|eq_0_main_thm_B_O},
the last inequality ensures that $C_{\phi}$ will be compact on $A_{\alpha}^{\psi}\left(\mathbb{B}_{N}\right)$
if, for any $B>0$, there exists a constant $A>0$ such that\begin{equation}
\psi^{-1}\left(\left(\psi\left(A\psi^{-1}\left(1/h^{N\left(\alpha\right)}\right)\right)\right)^{\frac{\alpha-\beta}{N\left(\alpha\right)}}.\psi\left(C_{\beta}\psi^{-1}\left(1/h^{N\left(\beta\right)}\right)\right)\right)\geq B\psi^{-1}\left(1/h^{N\left(\alpha\right)}\right),\label{eq|eq_3_main_thm_B_O}\end{equation}
for $h$ small emough. Putting $x=\psi^{-1}\left(1/h^{N\left(\alpha\right)}\right)$,
\prettyref{eq|eq_3_main_thm_B_O} is equivalent to\[
\psi\left(Bx\right)\leq\psi\left(Ax\right)^{\frac{\alpha-\beta}{N\left(\alpha\right)}}.\psi\left(C_{\beta}\psi^{-1}\left(\left(\psi\left(x\right)\right)^{N\left(\beta\right)/N\left(\alpha\right)}\right)\right),\]
which is in turn satisfied, using the convexity of $\psi$ and $C_{\beta}\geq1$,
if\[
\psi\left(Bx\right)^{N\left(\alpha\right)}\leq\psi\left(Ax\right)^{\alpha-\beta}.\psi\left(x\right)^{N\left(\beta\right)},\]
for $x$ large enough. Let us notice that this last inequality is
equivalent to\[
\frac{\psi\left(Bx\right)^{N\left(\beta\right)/\left(\alpha-\beta\right)}}{\psi\left(x\right)^{N\left(\beta\right)/\left(\alpha-\beta\right)}}\leq\frac{\psi\left(Ax\right)}{\psi\left(Bx\right)},\]
for $x$ large enough, which is nothing but the $\nabla_{0}$-Condition
(see Paragraph 2.2).\end{proof}
\begin{strem}
We mention that the proof of the necessary part of the previous theorem
does not use the boundedness of $C_{\phi}$ on some {}``smaller''
weighted Bergman-Orlicz space. Also, it is not necessary to assume
that $\psi$ satisfies the $\nabla_{0}$-Condition.
\end{strem}
Since every composition operator is bounded on every $A_{\alpha}^{\psi}\left(\mathbb{B}_{N}\right)$
as soon as $\psi$ satisfies the $\Delta^{2}$-Conditon, we have the
following corollary:
\begin{stcor}
Let $\alpha>-1$, let $\psi$ be an Orlicz function satisfying the
$\Delta^{2}$-Condition and let $\phi:\mathbb{B}_{N}\rightarrow\mathbb{B}_{N}$
be holomorphic. Then $C_{\phi}$ is compact on $A_{\alpha}^{\psi}\left(\mathbb{B}_{N}\right)$
if and only if\[
\lim_{\left|z\right|\rightarrow1}\frac{\psi^{-1}\left(1/\left(1-\left|\phi\left(z\right)\right|\right)\right)}{\psi^{-1}\left(1/\left(1-\left|z\right|\right)\right)}=0.\]
\end{stcor}
\begin{proof}
It is sufficient to remark that we have\[
\lim_{\left|z\right|\rightarrow1}\frac{\psi^{-1}\left(1/\left(1-\left|\phi\left(z\right)\right|\right)^{N\left(\alpha\right)}\right)}{\psi^{-1}\left(1/\left(1-\left|z\right|\right)^{N\left(\alpha\right)}\right)}=0\,\Longleftrightarrow\,\lim_{\left|z\right|\rightarrow1}\frac{\psi^{-1}\left(1/\left(1-\left|\phi\left(z\right)\right|\right)\right)}{\psi^{-1}\left(1/\left(1-\left|z\right|\right)\right)}=0,\]
since $\psi$ satisfies the $\Delta^{2}$-Condition. Indeed, it is
easy to deduce from the definition of the $\Delta^{2}$-Condition
that, if $a>1$, then $\psi^{-1}\left(x^{a}\right)\leq C\psi^{-1}\left(x\right)$
for some constant $C>0$ and for $x$ large enough.
\end{proof}
\smallskip{}

This corollary highlights an important difference with the classical
weighted Bergman case: when $\psi$ satisfies the $\Delta^{2}$-Condition,
the compactness (as well as the boundedness, \prettyref{thm|main_thm_B_O_rappel})
of composition operators on $A_{\alpha}^{\psi}\left(\mathbb{B}_{N}\right)$
does not depend on $\alpha>-1$.

Yet, this independency does not stand for $\alpha=-1$, i.e. for Hardy-Orlicz
spaces; indeed, it was shown (\cite[Theorem 5.8]{QUEF-LI-LE-RO-PI})
that there exists some Orlicz function $\psi$ which satisfies the
$\Delta^{2}$-Condition (to be precise, $\psi\left(x\right)=e^{x^{2}}-1$)
such that there exists a holomorphic self-map of $\mathbb{D}$ inducing
a compact operator on $A_{\alpha}^{\psi}\left(\mathbb{D}\right)$,
but not compact on $H^{\psi}\left(\mathbb{D}\right)$.

\smallskip{}

Nevertheless, the same proof as that of the necessary part of \prettyref{thm|thm_MAIN_3_B_O}
for Hardy-Orlicz spaces yields:
\begin{stprop}
Let $\phi:\mathbb{B}_{N}\rightarrow\mathbb{B}_{N}$ be holomorphic
and let $\psi$ be an Orlicz function. If $C_{\phi}$ is compact on
$H^{\psi}\left(\mathbb{B}_{N}\right)$, then\[
\lim_{\left|z\right|\rightarrow1}\frac{\psi^{-1}\left(1/\left(1-\left|\phi\left(z\right)\right|\right)\right)}{\psi^{-1}\left(1/\left(1-\left|z\right|\right)\right)}=0.\]

\end{stprop}
\bigskip{}

\end{document}